\documentclass[a4paper,11pt]{article}

\usepackage[utf8]{inputenc}
\usepackage[T1]{fontenc}
\usepackage{amsmath}
\usepackage{amsfonts}
\usepackage{amssymb}
\usepackage{amsthm}
\usepackage[version=4]{mhchem}
\usepackage{stmaryrd}
\usepackage{bbold}
\usepackage{graphicx}
\usepackage[export]{adjustbox}
\graphicspath{ {./images/} }

\usepackage{cite}

\usepackage{enumitem}
\usepackage{multibib}
\newcites{supp}{Supplementary References}

\newtheorem{theorem}{Theorem}
\newtheorem{corollary}[theorem]{Corollary}
\newtheorem{lemma}[theorem]{Lemma}
\newtheorem*{theorem*}{Theorem}
\newtheorem*{corollary*}{Corollary}
\newtheorem*{lemma*}{Lemma}
\def\bet{\begin{theorem}}
\def\eet{\end{theorem}}
\def\bel{\begin{lemma}}
\def\eel{\end{lemma}}
\def\bec{\begin{corollary}}
\def\eec{\end{corollary}}
\def\betn{\begin{theorem*}}
\def\eetn{\end{theorem*}}
\def\beln{\begin{lemma*}}
\def\eeln{\end{lemma*}}
\def\becn{\begin{corollary*}}
\def\eecn{\end{corollary*}}
\def\beq{\begin{equation}}
\def\eeq{\end{equation}}
\def\beqn{\begin{equation*}}
\def\eeqn{\end{equation*}}

\DeclareGraphicsExtensions{.eps,.jpg,.pdf,.ps} 

\title{Remarks on Derivations on Maximal Triangular Operator Algebras}

\author{Mark Spivack}
\date{\today}

\begin{document}

\maketitle

\section*{Abstract}

This note concerns bounded derivations on maximal triangular operator
algebras on a Hilbert space.   
Given any bounded derivation $\delta$ on a maximal triangular algebra
whose invariant lattice is continuous at 1, an operator which is shown
to implement $\delta$ 
is constructed explicitly. For a general reducible maximal triangular
algebra and any invariant projection $p$ the same  
construction yields an operator which implements $\delta$ on $pSp$
and $p^\perp S p^\perp$, and on all of $S$
if and only if $\delta$ obeys an additional triple product
rule and its restriction to $pSp^\perp$ is strongly continuous.

{\sl This work is based on unpublished parts of the author's
  dissertation \cite{spivack1982thesis} and describes a variant of the
  proof of a more general result in \cite{spivack1987certain} and is
  thus in effect a footnote to that work. The construction and the
  corresponding necessary and sufficient condition (strong continuity
  together with Property~$D_1$, of which the triple product rule here
  is a special case) were given in \cite{spivack1987certain}. To the
  best of the author's knowledge the result for algebras whose
  invariant lattice is continuous at 1 has not appeared elsewhere.
  The work here is not set in the context of the large body of
  subsequent research,   and no claims are made regarding its
  relationship to later developments.}

\section{Introduction}

Maximal triangular algebras were studied by Kadison and Singer
\cite{kadison2} and their properties have been the subject of numerous
subsequent studies. 
In this note we consider derivations on maximal 
triangular algebras and the questions of when, and by what operators,
they are automatically implemented. The results 
obtained are proofs by construction, i.e. an operator is exhibited
explicitly which implements the given derivation.

Derivations, on both bounded and unbounded operator algebras, have
been studied extensively (e.g. \cite{kadison1,
  ringrose2,christensen2,christensen3,bratteli2}), 
in particular on $C^*$ and von Neumann algebras
\cite{sakai3,johnson2,bratteli1}, 
and certain classes of non-self-adjoint algebras
\cite{christensen1,gilfeather}. 
Nest algebras, consisting of all operators in $B(H)$ leaving  
invariant every element of a given subspace nest and which thus share many
properties with maximal triangular algebras, were introduced by
Ringrose in \cite{ringrose1965some} (see also
\cite{lance,arveson,erdos,davidson,davidson1988nest,solel}). %
Christensen \cite{christensen1} showed that all
derivations on nest algebras are implemented. A constructive proof of
this result was given in \cite{spivack1982thesis,spivack1985derivations},
in which the cases of invariant lattices respectively continuous and discrete
at 1 were treated separately, and which motivates the approach here.

A maximal triangular algebra is {\sl reducible} if it has any
non-trivial invariant subspaces; otherwise it is irreducible.
For reducible algebras whose invariant lattices are strongly continuous
at $1$ we show by construction that any derivation is implemented.
Subsequently, for any reducible algebra $S$ with an invariant projection
$p$ we construct an operator which implements the given derivation on the
subsets $pSp$ and $p^{\perp}Sp^{\perp}$. Since $p^{\perp}Sp=\{0\}$, this
leaves only $pSp^{\perp}$, which equals $pB(H)p^{\perp}$ by Lemma~1. We
show that, for $\delta$ strongly continuous on $pSp^{\perp}$,
implementation of $\delta$ on this set is equivalent to a natural triple
product rule on $\delta$. This is a special case of the necessary and
sufficient condition of \cite{spivack1987certain} (strong continuity
together with  
Property~$D_1$).
For a non-trivial nest algebra the same construction gives an explicit
implementing operator, with no separate treatment of lattices continuous
and discrete at $1$; see Remark~2.

\section{Definitions and preliminaries}

Let $B(H)$ be the set of bounded linear operators on a complex Hilbert
space $H$.   An algebra $S \subset B(H)$ is {\sl triangular} if the algebra
${S} \cap {S}^{*}$ is maximal Abelian in ${B}({H})$. Then ${S} \cap
{S}^{*}$ is the {\sl diagonal} of ${S}$. 
For any triangular algebra $S_{1}$ containing $S$ we
have $S \cap S^{*}=S_{1} \cap S_{1}^{*}$ since $S \cap S^{*}$ is maximal Abelian.

Note that $\operatorname{lat}({S}) \subset {S} \cap {S}^{*}$ :  For $p
\in \operatorname{lat}(S)$ and $a \in S \cap S^{*}$, we have $a p=p a
p, a^* p=p a^* p$, so $p$ commutes with $S \cap S^{*}$. However $S \cap
S^{*}$ is maximal Abelian, so $p \in S \cap S^{*}$. For
any $p \in \operatorname{lat}(S), a \in S$, ~~$p^{\perp} a p=p^{\perp} p a
p=0$, so that ${p}^{\perp} {Sp}=\{0\}$.

A {\sl derivation on} an operator algebra $A$ is a
linear map $\delta$ from $A$ into $B(H)$ obeying the product rule
$\delta(ab)=\delta(a)b+a\delta(b)$.  
For a projection $p$ in the
domain of $\delta$, $\delta(p)$ is an operator which maps the range of
$p$ into its kernel, and vice versa, so that
$\left(\delta(p)\right)^2$ commutes with $p$.
(To verify this, consider the identity $\delta(p)=\delta(p^2)$ and
apply the product rule.)
We identify $p$ with
its range where this is unambiguous, and so for example may write
$\eta \in p$ when $p\eta=\eta$.

For any $b \in B(H)$ the map $d_{b}$ defined on a subalgebra $A$ of
$B(H)$ by
$d_{b}(a)=b a-a b$ for all $a \in A$ is a derivation. If $\delta$
is a derivation on $A$ and $\delta=d_{b}$ for some $ {b} \in B(H)$ then we
say that $\delta$ is {\sl implemented}. If also
$\delta: {A} \rightarrow {A}$ then $\delta$ is {\sl inner}.

If $\xi, \eta$ are non-zero vectors in $H$, then $<\xi, \eta>$
denotes the inner-product of $\xi$ with $\eta$, and $\xi \otimes \eta$
denotes the rank one operator given by
$(\xi \otimes \eta) \zeta = <\zeta, \xi> \eta$ for all $\zeta$ in
$ {H}$.

We require a few known properties and results (for most of which we
omit proofs):  

\bel\label{lemma1}
\cite{kadison2}. If ${p} \in
\operatorname{lat}({S})$ then ${pB}({H}) {p}^{\perp} \subset {S}$.  
\eel

\bec\label{cor2}
\cite{kadison2} If ${p} \in \operatorname{lat}({s})$ then ${S} \xi
\supset {p} ~~\forall \xi \in {p}^{\perp}$. 
\eec

\bel\label{lemma3}
\cite{kadison2} Lat(S) is totally-ordered.
\eel

\bel
The commutant $S'$ of $S$ is trivial i.e. $S'=\mathbb{C}$.
\eel
\noindent Hence, from a trivial calculation,  for any derivation
$\delta$ on $S$, any two linear  
maps which implement $\delta$ must differ by some $\alpha \in \mathbb{C}$.  

We will also need the following result for derivations acting on
projections:
\bel\label{2.3.1}
\cite{spivack1982thesis}
Let $ {p}$ be any projection and let $\delta$ be any
derivation whose domain includes $ {p}$. 
Then \\ (i) $\delta(p)$ is an operator mapping $r(p)$ into $n(p)$ and
vice versa, \\
(ii) $\delta$ is implemented on $p$ by the operator
$
b=(1-2 p) \delta(p) .
$
and \\ (iii) If $\delta_{1}(p)$ denotes any operator which maps $r(p)$
into $n(p)$ and vice versa then the linear extension of $\delta_{1}$
to the algebra (p) generated by $p$ and $1$ is a derivation. 
\eel
\noindent{\bf Proof.} (i) $p^{2}=p$, so that
\beqn
\delta(p)=p \delta(p)+\delta(p) p .
\eeqn
Pre-multiplying by $p$ we get $p \delta(p) p=0$, and since $\delta$ is
defined also on $p^{\perp}=1-p$ we get a similar result for
$p^{\perp}$, so that $ {p} \delta( {p})  {p}= {p}^{\perp} \delta( {p})
{p}^{\perp}=0$. 

(ii) Compute, using (i):
$(1-2 p) \delta(p) p-p(1-2 p) \delta(p)=\delta(p) p+p \delta(p)=\delta(p)$.

(iii) (p) is just the set $\{\alpha+\beta p: \alpha, \beta \in
\mathbb{C}\}$, and since $\delta_{1}(\mathbb{C})$ must be zero,
$\delta_{1}$ is defined on $\alpha+\beta p \in(p)$ by
$\delta_{1}(\alpha+\beta p)=\beta \delta_{1}(p)$. Also
$\delta_{1}$ is automatically linear.

Let $\alpha_{1}, \alpha_{2}, \beta_{1}, \beta_{2} \in \mathbb{C}$.
Then $\delta_{1}\left[\left(\alpha_{1}+\beta_{1}
    p\right)\left(\alpha_{2}+\beta_{2}
    p\right)\right]=\delta_{1}\left[\alpha_{1}
  \alpha_{2}+\left(\alpha_{1} \beta_{2}+\alpha_{2} \beta_{1}+\beta_{1}
    \beta_{2}\right) p\right]$ 
$=\left(\alpha_{1} \beta_{2}+\alpha_{2} \beta_{1}+\beta_{1}
  \beta_{2}\right) \delta_{1}(p)$. 
and $\delta_{1}\left(\alpha_{1}+\beta_{1}
  p\right)\left(\alpha_{2}+\beta_{2}
  p\right)+\left(\alpha_{1}+\beta_{1} p\right)
\delta_{1}\left(\alpha_{2}+\beta_{2} p\right)$ 

$=\alpha_{2} \beta_{1} \delta_{1}(p)+\beta_{1} \beta_{2} \delta_{1}(p)
p+\alpha_{1} \beta_{2} \delta_{1}(p)+\beta_{1} \beta_{2} p
\delta_{1}(p)$ 

$=\left(\alpha_{1} \beta_{2}+\alpha_{2} \beta_{1}+\beta_{1}
  \beta_{2}\right) \delta_{1}(p)$. 
Hence $\delta_{1}$ is a derivation.   \qed

\section{Derivations on maximal triangular algebras}\label{sec4.3}

Throughout the remainder $S$ will denote a reducible maximal triangular algebra.
Let $\delta: S \rightarrow B(H)$ be a continuous derivation.
For the case in which the invariant lattice is strongly continuous at
1 an operator is constructed which implements 
$\delta$.  We also formulate a condition on $\delta$
under which, for any reducible maximal triangular algebra,
the constructed operator is well-defined and implements $\delta$.  

We now proceed with the construction. There is $p
\in \operatorname{lat}(S)$ such that $p \neq 0$ or 1 . Let $\xi \in
p^{\perp}, \eta \in p$ be unit vectors. Then by Lemma \ref{lemma1} the rank one
operator $\xi \otimes \eta$ is in $S$ and $S \xi \supset p$. 
This immediately gives the first result:

\bel\label{lemma9} 
For  $p \in \operatorname{lat}(S)$ such that $p \neq 0$ or 1, 
choose unit vector $\xi_{0} \in {p}^{\perp}$ and write
${p}_{0}=\xi_{0} \otimes \xi_{0}$. 
Define a map $b_1$ by 
\beq
b_1 {a} \xi_{0}=\delta\left({ap}_{0}\right) \xi_{0}
\label{implement1}
\eeq
where $a \in {S}$, $a \xi_{0} \in {p} $, and $b_1 {p}^{\perp}=0$.
Then ${b_1}$
is a well-defined and bounded linear operator, with 
$\left\|b_1 \right\| \leqslant\|\delta\|$, and $\delta({a})
{p}={d}_{b_1}(a){p} ~~\forall {a} \in {S}$. 
\eel

\noindent{\bf Proof.} 
Although ${p}_{0}$ is not
necessarily contained in ${S}$, if a $\in {S}$ such that $a \xi_{0}
\in {p}$ then $a p_{0}=p\left(a p_{0}\right) p^{\perp} \in S$ and
$\left\|a p_{0}\right\|=\left\|a \xi_{0}\right\|$. Also $a
\xi_{0}=0 \Rightarrow a p_{0}=0$. 
Thus it follows immediately that ${b_1}$ is well-defined and linear.
Furthermore
\beqn
\left\|{b_1} {a} \xi_{0}\right\|=\left\|\delta\left({ap}_{0}\right)
  \xi_{0}\right\|
\leqslant\|\delta\|\left\|{ap}_{0}\right\|=\|\delta\|\left\|{a}
  \xi_{0}\right\| 
\eeqn
whenever $a \xi_{0} \in {p}$, and so ${b_1} \in {B}({H})$ and
$\left\|{b_1}\right\| \leqslant\|\delta\|$. 
Let $a, c \in S$ such that $c \xi_{0} \in p$.
Then $ac \xi_{0} \in {p}$, ${cp}_{0} \in {S}$ and
\beqn
\left({b_1} a-a {b_1}\right) c \xi_{0}  =\left[\delta\left(a c
    p_{0}\right)-a \delta\left(c p_{0}\right)\right] \xi_{0}  
 =\delta(a) c \xi_{0} .
\eeqn
Hence ${b_1}$ implements $\delta$ on ${p}$. \qed

We can now state the main derivation result for those algebras $S$
with the property that lat($S$) is strongly continuous at 1. 

\bet\label{theorem10}
Suppose ${S}$ is such that 1 is the strong limit of
projections in lat($S$). Then $\delta$ is implemented and we can
construct $b \in B(H)$ such that 
$\delta=d_{b}\mid_ S$ and $\|b\| \leqslant 2\|\delta\|$. 
\eet
\noindent{\bf Proof.}
Choose $p \in \operatorname{lat}(S)$ with $p \neq 0$. Since 1 is the strong
limit of projections in $\operatorname{lat}(S)$, there is $p_1 \in
\operatorname{lat}(S)$ with $p < p_1 < 1$; choose a unit vector $u \in
p_1 \ominus p$. Let $b_0$ be the operator given by Lemma \ref{lemma9}
for $p$, constructed with $\xi_0 = u$. For each $p_\alpha \in
\operatorname{lat}(S)$ with $p_1 \leqslant p_\alpha < 1$ choose a unit
vector $v_\alpha \in p_\alpha^{\perp}$, and let $b_\alpha$ be the
operator given by Lemma \ref{lemma9} for $p_\alpha$, constructed with
$\xi_0 = v_\alpha$, so that $\|b_\alpha\| \leqslant \|\delta\|$ and
$\delta(a) p_\alpha = d_{b_\alpha}(a) p_\alpha$ for all $a \in S$.

Let $\eta \in p$ be a unit vector. Taking $a = \xi_0 \otimes \eta$ in
(\ref{implement1}) gives
\beqn
b_0 \eta = \delta(u \otimes \eta) u , \qquad
b_\alpha \eta = \delta(v_\alpha \otimes \eta) v_\alpha .
\eeqn
By Lemma \ref{lemma1}, $u \otimes \eta \in S$ (since $u \in p^{\perp}$,
$\eta \in p$) and $v_\alpha \otimes u \in S$ (since $v_\alpha \in
p_\alpha^{\perp}$, $u \in p_\alpha$), and
$(u \otimes \eta)(v_\alpha \otimes u) = v_\alpha \otimes \eta$. Applying
the product rule and evaluating at $v_\alpha$,
\beqn
\delta(v_\alpha \otimes \eta) v_\alpha
= \delta(u \otimes \eta) u + (u \otimes \eta)\, \delta(v_\alpha \otimes u)
v_\alpha
= \delta(u \otimes \eta) u + \mu_\alpha \eta ,
\eeqn
where $\mu_\alpha = < \delta(v_\alpha \otimes u) v_\alpha , u >$ does
not depend on $\eta$. Hence $(b_\alpha - b_0) p = \mu_\alpha p$, and
$|\mu_\alpha| \leqslant \|\delta\|$ since $\|v_\alpha \otimes u\| = 1$.

Put $b_\alpha^{\prime} = b_\alpha - \mu_\alpha p_\alpha$. Since $a
p_\alpha = p_\alpha a p_\alpha$ for $a \in S$, $b_\alpha^{\prime}$ still
implements $\delta$ on $p_\alpha$. Moreover $b_\alpha^{\prime} p = b_0
p$ and
\beqn
\|b_\alpha^{\prime}\| \leqslant \|b_\alpha\| + |\mu_\alpha|
\leqslant 2 \|\delta\| .
\eeqn

Now let $p_1 \leqslant p_\alpha < p_\beta < 1$ and put $T =
(b_\alpha^{\prime} - b_\beta^{\prime}) p_\alpha$. Since both operators
implement $\delta$ on $p_\alpha$, $T a x = a T x$ for all $a \in S$ and
$x \in p_\alpha$. Taking $a = p_\alpha$ and $a = p$ shows that $T$ maps
$p_\alpha$ into itself and commutes with $p$ there. Taking $a = \xi
\otimes \eta$ with $\xi \in p_\alpha \ominus p$ and $\eta \in p$ unit
vectors (so $a \in S$ by Lemma \ref{lemma1}) gives
$< x , \xi > T \eta = < T x , \xi > \eta$ for all $x \in p_\alpha$. With
$x = \xi$ this shows that $T$ is a scalar $\lambda$ on $p$, and then
$< T x , \xi > = \lambda < x , \xi >$ for all such $\xi$, so that $T =
\lambda$ on $p_\alpha \ominus p$ as well. Thus $T = \lambda p_\alpha$.
But $T p = (b_0 - b_0) p = 0$, so $\lambda = 0$, that is
$b_\alpha^{\prime} p_\alpha = b_\beta^{\prime} p_\alpha$.

We may therefore define $b$ on the dense subspace $\bigcup_\alpha
p_\alpha H$ by $b x = b_\alpha^{\prime} x$ for $x \in p_\alpha$, and
since $\|b_\alpha^{\prime}\| \leqslant 2 \|\delta\|$ for all $\alpha$,
$b$ extends to an operator in $B(H)$ with $\|b\| \leqslant 2
\|\delta\|$; equivalently $b$ is the strong limit of
$\{b_\alpha^{\prime}\}$. For $a \in S$ we have $a p_\alpha = p_\alpha a
p_\alpha$, so
\beqn
d_b(a) p_\alpha = b_\alpha^{\prime} a p_\alpha - a b_\alpha^{\prime}
p_\alpha = \delta(a) p_\alpha ,
\eeqn
and letting $p_\alpha \rightarrow 1$ strongly gives $\delta = d_b \mid
S$. \qed

For $p \in \operatorname{lat}(S)$, $p \ne 0$ or $1$, Lemma \ref{lemma9} gives an
operator $b_1 = b_1p$. Continuing with the construction we define
$c_{1}=p c_{1} p^{\perp}$ 
by
\beq\label{c_1}
c_{1} p^{\perp}=-\delta(p) p^{\perp}, ~~~~~
c_{1} p=0 ~~.
\eeq
Then ${pc}_{1} {p}^{\perp}={c}_{1} {p}^{\perp}$ by Lemma \ref{2.3.1}. It
is clear that ${c}_{1} \in {B}({H})$ and $\left\|c_{1}\right\|
\leqslant\|\delta\|$. Now define $b_{2}$ by
\beq
b_{2}={b_1}+c_{1} .
\label{b_2}
\eeq
Then we have

\bet\label{theorem11}
With $c_1$, and $b_2$ defined as in (\ref{c_1}),(\ref{b_2}),   $b_{2}$
implements $\delta$ 
on the algebra 
${Sp}={pSp}$, and $\left\|{b}_{2}\right\| \leqslant 2\|\delta\|$. 
\eet
\noindent{\bf Proof.} Let a $\in$ S. With the notation of the construction of
${b_1}$, if $\xi \in {p}$ then 
\beqn
\begin{aligned}
\left(b_{2} a p-a p b_{2}\right) \xi & =\left(b_{2} p a p-a p
  b_{2}\right) \xi \\ 
& =\delta\left(a p c p_{0}\right) \xi_{0}-a p \delta\left(c
  p_{0}\right) \xi_{0} ,~~~~~~\text { where }~ c \in S, ~~c
\xi_{0}=\xi, \\  
& =\delta(a p) c p_{0} \xi_{0}=\delta(a p) \xi .
\end{aligned}
\eeqn
If $\xi \in p^{\perp}$ then
\beqn
\begin{aligned}
\left({b}_{2} {ap}-{ap} {b}_{2}\right) \xi & ~~=~~-{apb}_{2}
\xi~~=~~{ap} \delta({p}) {p}^{\perp} \xi \\ 
& ~~=~~{ap} \delta({p}) \xi~~=~~\delta({ap}) \xi-\delta({ap}) {p}
\xi~~=~~\delta({ap}) \xi . 
\end{aligned}
\eeqn
Hence $\delta={d}_{{b}_{2}} \mid {Sp}$.  \qed

\noindent{\bf Remark.} ${b}_{2}$ clearly implements $\delta$ on ${p}$ since
\beqn
\left(b_{2} a-a b_{2}\right) p  = b_{2} p a p-a b_{2} p 
 = \left({b_1} a-a {b_1}\right) p \quad \forall a \in S,
\eeqn
so we could have stated Lemma \ref{lemma9} with $b_{2}$ instead of ${b_1}$.

We come now to the final part of the construction.
Recall that if $q$ is the rank one operator $\xi \otimes \eta$ then
$q^{*}$ is the rank one operator $\eta \otimes \xi$.

If $p_{1}$ is the rank one projection $\xi \otimes \xi$ then
$p_{1}=q^{*} q$ and, in general, 

\beqn
\left(\eta \otimes \xi_{2}\right)\left(\xi_{1} \otimes \eta\right)=\xi_{1}
\otimes \xi_{2} . 
\eeqn

This brings us to the following: 
\bel\label{lemma6_c2}
Choose a fixed unit vector $\eta_{1} \in p$ and put $q_{1}=\xi_{0}
\otimes \eta_{1}$ so that $q_{1} \in S$. 
Define a map $c_2$ as follows:
\beq
c_{2} p = 0 \text {, and } p c_{2} p^{\perp}=0
\eeq
and for any $\xi \in p^{\perp}$
\beq
p_{\xi} c_{2} p^{\perp}  =-q^* \delta(q) p^{\perp}+q^*
\delta\left(q_{1}\right) {q_{1}}^* q
\eeq
where $q=\xi \otimes \eta_{1} = q p^{\perp}$, and ${p}_{\xi}=\xi \otimes \xi$.
Then  ${c}_{2} \in {B}({H})$.
\eel
(Notice that this definition depends on the fixed vectors $\xi_{0} \in
{p}^{\perp}$, $\eta_{1} \in$ p. It is not obvious that the definition is
invariant to within a constant additive factor under these choices.)

\noindent{\bf Proof.}
${c}_{2}$ is clearly well-defined. It is easy to
see that $c_{2}$ is bounded and that for any unit vector
$\eta,\left\|c_{2} \eta\right\| \leqslant 2\|\delta\|$ since ${pc}_{2}=0$
and, for any unit vector $\xi \in {p}^{\perp},\left\|p_{\xi} c_2
  \eta\right\| \leqslant 2\|\delta\|$.  

It remains to show that ${c}_{2}$ is linear. For this it must be shown
that for unit vectors $\xi,\left\{\xi_{i}\right\}_{i \in I}$ in
$p^{\perp}$ such that $\xi=\sum_{i} \alpha_{i} \xi_{i}$, say, the
definitions ${p}_{\xi_{i}} {c}_{2} {p}^{\perp}$ as above lead to the
same value for $p_{\xi} c_{2} {p}^{\perp}$ as by defining this
directly. It will suffice to show this for, say, $\xi=r_{\alpha}
\xi_{\alpha}+r_{\beta} \xi_{\beta}$,  where
$\|\xi\|=\left\|\xi_{\alpha}\right\|=\left\|\xi_{\beta}\right\|=1$ and
$r_{\alpha}$, $r_{\beta} \in \mathbb{C}$. 

Write ${p}_{\alpha}=\xi_{\alpha}{\otimes \xi_{\alpha}}, ~~
{p}_{\beta}=\xi_{\beta}{\otimes \xi_{\beta}},~~ {p}_{\xi}=\xi \otimes
\xi, ~~q_{\alpha}=\xi_{\alpha}{\otimes \eta_{1}},~~
q_{\beta}=\xi_{\beta}{\otimes \eta_{1}}, ~~ \text{and}~~ q=\xi \otimes
\eta_{1} .$

We have first that, $\forall \zeta \in H$,

\beqn
q \zeta= < \zeta, r_{\alpha} \xi_{\alpha}+r_{\beta}
  \xi_{\beta}>\eta=\left(\bar{r}_{\alpha} q_{\alpha}+\bar{r}_{\beta}
    q_{\beta}\right) \zeta . 
\eeqn

This gives $q^{*}=r_{\alpha} q_{\alpha}^{*}+r_{\beta} q_{\beta}^{*}$ and

\beqn
p_{\xi}=q^{*} q=\left|r_{\alpha}\right|^{2}
p_{\alpha}+\left|r_{\beta}\right|^{2} p_{\beta}+r{ }_{\alpha}
\bar{r}_{\beta} q_{\alpha}^{*} q_{\beta}+\bar{r}_{\alpha} r_{\beta}
q_{\beta}^{*} q_{\alpha} \cdot 
\eeqn

Defining ${p}_{\xi} {c}_{2} {p}^{\perp}$ directly we thus have
\begin{align}
p_{\xi} c_{2} p^{\perp} & =q^{*} \delta\left(q_{1}\right) q_{1}^{*}
                          q-q^{*} \delta(q) p^{\perp} \\ 
& =\left|r_{\alpha}\right|^{2} q_{\alpha}^{*} \delta\left(q_{1}\right)
  q_{1}^{*}   q_{\alpha} + \left|r_{\beta}\right|^{2}  q_{\beta}^{*} 
 \delta\left(q_{1}\right) q_{1}^{*} q_{\beta} + r_{\alpha} \bar{r}_{\beta}
  q_{\alpha}^{*}  \delta(q_{1}) q_{1}^{*} q_{\beta}
  \\  
& +\bar{r}_{\alpha} r_{\beta} q_{\beta}^{*} \delta\left(q_{1}\right)
         q_{1}^{*} q_{\alpha} \\ 
& -\left[\left|r_{\alpha}\right|^{2} q_{\alpha}^{*} \delta\left(q_{\alpha}\right)
 +\left|r_{\beta}\right|^{2} q_{\beta}^{*} \delta\left(q_{\beta}\right)+r_{\alpha} \bar{r}_{\beta} q_{\alpha}^{*} \delta\left(q_{\beta}\right)+\bar{r}_{\alpha} r_{\beta} q_{\beta}^{*} \delta\left(q_{\alpha}\right)\right] \label{1}
\end{align}

On the other hand

\beqn
\begin{aligned}
& {p}_{\xi} {c}_{2} {p}^{\perp}=\left[\left|r_{\alpha}\right|^{2}
  {p}_{\alpha}+\left|r_{\beta}\right|^{2} {p}_{\beta}+r_{\alpha}
  \bar{r}_{\beta} {q}_{\alpha}^{*} {q}_{\beta}+\bar{r}_{\alpha}
  r_{\beta}q_{\beta}^{*} {q}_{\alpha}\right] b p^{\perp} \\ 
& =\left|r_{\alpha}\right|^{2}\left[q_{\alpha}^{*}
  \delta\left(q_{1}\right) q_{1}^{*} q_{\alpha}-q_{\alpha}^{*}
  \delta\left(q_{\alpha}\right)
  p^{\perp}\right]+\left|r_{\beta}\right|^{2}\left[q_{\beta}^{*}
  \delta\left(q_{1}\right) q_{1}^{*} q_{\beta}-q_{\beta}^{*}
  \delta\left(q_{\beta}\right) p^{\perp}\right] 
\end{aligned} 
\eeqn

\begin{align}
& +r_{\alpha} \bar{r}_{\beta} q_{\alpha}^{*} q_{\beta}\left[q_{\beta}^{*} \delta\left(q_{1}\right) q_{1}^{*} q_{\beta}-q_{\beta}^{*} \delta\left(q_{\beta}\right) p^{\perp}\right] \\
& +\bar{r}_{\alpha} r_{\beta} q_{\beta}^{*} q_{\alpha}\left[q_{\alpha}^{*} \delta\left(q_{1}\right) q_{1}^{*} q_{\alpha}-q_{\alpha}^{*} \delta\left(q_{\alpha}\right) p^{\perp}\right] \label{2}
\end{align}

Notice that $q_{\alpha}^{*} q_{\beta} q_{\beta}^{*}=q_{\alpha}^{*}$
and $q_{\beta}^{*} q_{\alpha} q_{\alpha}^{*}=q_{\beta}^{*}$, and these
occur in the last two terms of the expression above. Comparing (\ref{1}) and
(\ref{2}), then, term by term we see that they are equal. 
 
This extends by induction to any finite sum of vectors in
${p}^{\perp}$ and hence, since we have shown $c_{2}$ to be bounded, to
any sum in $p^{\perp}$. It follows that ${c}_{2}$ is linear and the
result follows.   
\qed

We now have

\bet\label{theorem12}
Define $b=b_{2}+c_{2}$, with $b_2$ as in equation (\ref{b_2}) and
$c_2$ defined as in
Lemma \ref{lemma6_c2}. 
Then b implements $\delta$ on the algebras ${Sp}$ and
${p}^{\perp} {Sp}^{\perp}$, and 
$\|b\| \leqslant 4\|\delta\|$.
\eet
\noindent{\bf Proof.} We already have that
$\left\|b_{2}\right\|,\left\|c_{2}\right\| \leqslant 2\|\delta\|$, so
$\|b\| \leqslant 4\|\delta\|$. 

Let $a \in S$.
Then $(bpap-papb) ~=~ {b}_{2} ~ pap-papb_{2}=\delta( pap )$ by
Theorem \ref{theorem11}. 
To show that $\delta={d}_{{b}} \mid {p}^{\perp} {Sp}^{\perp}$ we
consider various cases. 

\begin{enumerate}[label=\arabic*)]
\item Let $\xi \in p$.
Then $\left(b p^{\perp} a p^{\perp}-p^{\perp} a p^{\perp} b\right)
\xi=-p^{\perp} a p^{\perp} b \xi=-p^{\perp} a p^{\perp} {b_1} \xi$ 
\beqn
\begin{aligned}
& =-p^{\perp} a p^{\perp} \delta(q) \xi_{0}, ~~~~~\text { where }
q=\xi_{0} \otimes \xi \\ 
& =\delta\left(p^{\perp} a p^{\perp}\right) q
\xi_{0}-\delta\left(p^{\perp} a p^{\perp} q\right)
\xi_{0}=\delta\left(p^{\perp} a p^{\perp}\right) \xi . 
\end{aligned}
\eeqn

 \item  Let $\xi \in p^{\perp}$ and consider
 $\operatorname{pd}_{{b}}\left({p}^{\perp} {ap}{ }^{\perp}\right)$ : 
\beqn
\begin{aligned}
& {p}\left(b {p}^{\perp} {a} {p}^{\perp} \xi-{p}^{\perp} {a}
  {p}^{\perp} {b} \xi\right)={pb} {p}^{\perp} {a} {p}^{\perp}
\xi={b}_{2} {p}^{\perp} {ap^{ \perp } \xi} \\ 
& =-\delta(p) p^{\perp} a p^{\perp} \xi=\delta\left(p^{\perp}\right)
p^{\perp} a p^{\perp} \xi \\ 
& =\delta\left(p^{\perp} a p^{\perp}\right) \xi-p^{\perp}
\delta\left(p^{\perp} a p^{\perp}\right) \xi=p \delta\left(p^{\perp} a
  p^{\perp}\right) \xi . 
\end{aligned}
\eeqn

\item Finally let $\xi \in p^{\perp}$ and consider $p^{\perp}
d_{b}\left(p^{\perp} a p^{\perp}\right)$ : \\
Choose a basis $\left\{\xi_{\alpha}\right\}$ for $p^{\perp}$ and put
$p_{\alpha}=\xi_{\alpha} \otimes \xi_{\alpha}$, 
$q_{\alpha}=\xi_{\alpha} \otimes \eta_{1} ~~~\forall \alpha$. 
We consider, for each $\alpha$, ~~$p_{\alpha} d_{b}\left(p^{\perp} a
  p^{\perp}\right)$. 
First we define the constants $r_{\alpha \beta} \in \mathbb{C}$ by
\beqn
p_{\alpha} a p_{\beta}=r_{\alpha \beta} q_{\alpha}^{*} q_{\beta} .
\eeqn
(It is easily seen that ${p}_{\alpha} {ap}_{\beta}$ and
${q}_{\alpha}^{*} {q}_{\beta}$ are multiples of the rank one operator
$\xi_{\beta} \otimes \xi_{\alpha}$.)
Note that
$ q_{\alpha} a p_{\beta}=q_{\alpha} p_{\alpha} a p_{\beta}=r_{\alpha \beta} q_{\beta} 
$, 
and similarly ${p}_{\alpha} {aq}_{\beta}^{*}={p}_{\alpha} {ap}_{\beta}
{q}_{\beta}^{*}={r}_{\alpha \beta} {q}_{\alpha}^{*}$.  

Then  ${p}_{\alpha}\left({bp}{ }^{\perp} {ap}{ }^{\perp}
  \xi-{p}^{\perp} {ap}{ }^{\perp} {b} \xi\right)$ 
\beqn
\begin{aligned}
& ={q}_{\alpha}^{*} \delta\left({q}_{1}\right) {q}_{1}^{*} {q}_{\alpha} {ap}{ }^{\perp} \xi-{q}_{\alpha}^{*} \delta\left({q}_{\alpha}\right) {p}^{\perp} {ap}{ }^{\perp} \xi-\sum_{\beta} {p}_{\alpha} {p}^{\perp} {a} {p}_{\beta} {b} \xi \\
& =-q_{\alpha}^{*} \delta\left(q_{\alpha}\right) p^{\perp} a p^{\perp} \xi+\sum_{\beta}\left[q_{\alpha}^{*} \delta\left(q_{1}\right) q_{1}^{*} q_{\alpha} a p_{\beta} \xi\right. \\
& ~~~~~~\left.-{p}_{\alpha} {aq}{ }_{\beta}^{\star} \delta\left({q}_{1}\right) {q}_{1}^{\star} {q}_{\beta} \xi+{p}_{\alpha} {a} {q}_{\beta}^{\star} \delta\left({q}_{\beta}\right) \xi\right] \\
& =-q_{\alpha}^{*} \delta\left(q_{\alpha}\right) p^{\perp} a p^{\perp} \xi \\
& ~~~~~~ +\sum_{\beta}\left[r_{\alpha \beta} q_{\alpha}^{*} \delta\left(q_{1}\right) q_{1}^{*} q_{\beta} \xi-r_{\alpha \beta} q_{\alpha}^{*} \delta\left(q_{1}\right) q_{1}^{*} q_{\beta} \xi+r_{\alpha \beta} q_{\alpha}^{*} \delta\left(q_{\beta}\right) \xi\right] \\
& =-q_{\alpha}^{*} \delta\left(q_{\alpha}\right) p^{\perp} a p^{\perp} \xi+\sum_{\beta} q_{\alpha}^{*} \delta\left(r_{\alpha \beta} q_{\beta}\right) \xi \\
& =-q_{\alpha}^{*} \delta\left(q_{\alpha}\right) p^{\perp} a p^{\perp} \xi+\sum_{\beta} q_{\alpha}^{*} \delta\left(q_{\alpha} a p_{\beta}\right) \xi \\
& =-q_{\alpha}^{\star} \delta\left(q_{\alpha}\right) p^{\perp} a p^{\perp} \xi+q_{\alpha}^{\star} \delta\left(q_{\alpha} p^{\perp} a p^{\perp}\right) \xi \\
& =q_{\alpha}^{*} q_{\alpha} \delta\left(p^{\perp} a p^{\perp}\right) \xi=p_{\alpha} \delta\left(p^{\perp} a p^{\perp}\right) \xi .
\end{aligned}
\eeqn
\end{enumerate}
Hence $\forall \alpha$, $p_{\alpha} d_{b}\left(p^{\perp} a
  p^{\perp}\right) \xi=p_{\alpha} \delta\left(p^{\perp} a
  p^{\perp}\right) \xi$, and the proof is complete. 
\qed

Since ${p}$ is invariant under ${S}$, and since ${p}^{\perp}
Sp=\{0\}$ and $S p=p S p$,  we can write $S=S  
p+pSp^{\perp}+p^{\perp} S p^{\perp}$. Theorem \ref{theorem12}
shows that $\delta$ is
implemented by ${b}$ on ${Sp}$ and ${p}^{\perp} {Sp}^{\perp}$,
so extending the result to $S$ itself hinges on whether it holds for
the remaining term ${pSp}^{\perp}$.  This we have not proved in general.
However, we give as the final result a natural condition on
$\delta$ in the form of triple product rule under which it does indeed hold.

To motivate this consider operators ${s}_{1}, {s}_{2}, {s}_{3} \in
{S}$ and ${a} \notin {S}$ such that $s_2 a s_3 \in {S}$. For
example if $s_{3}={s}_{3} {p}^{\perp}$ and ${s}_{2}={p} {s}_{2}$ then
by Lemma \ref{lemma1} $a$ can be any operator in $B(H)$. Then although
$\delta(a)$ 
is undefined one may ask whether the equation 
\beq
\delta\left(s_{1} s_{2} a s_{3}\right)=s_{1} \delta\left(s_{2} a
  s_{3}\right)+\delta\left(s_{1}\right) s_{2} a s_{3}  
\label{condition}
\eeq
holds. In general this would seem to be a very strong restriction, although
of course whenever $\delta$ is implemented it must hold since
an implemented derivation extends automatically to $B(H)$. 
The result concerns a less restrictive form of
condition (\ref{condition}).

\bet[Triple product rule]
\label{theorem13}

Suppose the vectors $\xi_{0} \in p^{\perp}$ and $\eta_{1} \in p$, the
basis $\left\{\xi_{\alpha}\right\}$ for $p^{\perp}$, and the operator
$b$ are all as previously defined. Let $\eta$ be any vector in ${p}$. 

Define the rank one operators
\beq
q_{\alpha}=\xi_{\alpha} \otimes \eta_{1},~~ q_{1}=\xi_{0} \otimes \eta_{1},
~~q=\xi_{\alpha} \otimes \eta \text {, and } ~~q_{2}=\xi_{0} \otimes \eta \cdot 
\label{q_defn}
\eeq

Then $\delta$ is implemented if and only if its restriction to
$pS{p}^{\perp}$ is strongly continuous and  
\begin{equation}
\delta(q)=\delta\left(q q_{\alpha}^{*} q_{1}\right) q_{1}^{*}
q_{\alpha}+q q_{\alpha}^{*} \delta\left(q_{\alpha}\right)-q
q_{\alpha}^{*} \delta\left(q_{1}\right) q_{1}^{*} q_{\alpha} .
\label{theorem13eq} 
\end{equation}
If so then $\delta=d_{b} \mid S$.
\eet
\noindent Here, as in \cite{spivack1987certain}, strong continuity means
continuity with respect to the strong operator topology on bounded
subsets of $pSp^{\perp}$.

\noindent{\bf Proof.}
We have $q=q q_{\alpha}^{*} q_{1} q_{1}^{*} q_{\alpha}$ and
$q_{2}=q q_{\alpha}^{*} q_{1}$. 
All the operators $q$, $q_\alpha$, $q_1$, $q_2$ and none of their
adjoints, are in ${S}$. 

In addition to the expressions for $q, q_{2}$ in
terms of the 
other operators we have $q_{\alpha}=q_{\alpha} q_{1}^{*} q_{1}$
and other such relations which arise wherever the "end-points" of the
various operators coincide. 
So (\ref{theorem13eq}) can be written
\beqn
\delta\left(q q_{\alpha}^{*} q_{1} q_{1}^{*}
  q_{\alpha}\right)=\delta\left(q q_{\alpha}^{*} q_{1}\right)
q_{1}^{*} q_{\alpha}+q q_{\alpha}^{*} \delta\left(q_{1} q_{1}^{*}
  q_{\alpha}\right)  
-q q_{\alpha}^{*} \delta\left(q_{1}\right) q_{1}^{*} q_{\alpha} \cdot
\eeqn
Assuming that $\delta$ is implemented, say by $c$, $\delta=d_{c}$ can be
defined on the whole of ${B}({H})$, and we can compute using $\delta$
as if this is so: 
 
\beqn
\begin{aligned}
\delta\left(q q_{\alpha}^{*} q_{1} q_{1}^{*} q_{\alpha}\right) &
=\delta\left(q q_{\alpha}^{*} q_{1}\right) q_{1}^{*} q_{\alpha}+q
q_{\alpha}^{*} q_{1} \delta\left(q_{1}^{*} q_{\alpha}\right) \\ 
& =\delta\left(q q_{\alpha}^{*} q_{1}\right) q_{1}^{*} q_{\alpha}+q
q_{\alpha}^{*} \delta\left(q_{1} q_{1}^{*} q_{\alpha}\right)-q
q_{\alpha}^{*} \delta\left(q_{1}\right) q_{1}^{*} q_{\alpha} . 
\end{aligned}
\eeqn
Moreover $x \mapsto cx-xc$ is strongly continuous on bounded subsets
of $B(H)$, so the restriction of $\delta$ to $pSp^{\perp}$ is strongly
continuous.

Suppose now that (\ref{theorem13eq}) holds.

If $a=\operatorname{pap}^{\perp} \in {S}$, then ${a}=\sum_{\alpha}
{pap}_{\alpha}$  where ${p}_{\alpha}=\xi_{\alpha} \otimes
\xi_{\alpha} \forall \alpha$ and each
${pap}_{\alpha} \in S$, the sum converging strongly with partial
sums bounded in norm by $\|a\|$.
Since $\delta$ and $d_b$ are both linear and strongly continuous on
bounded subsets of ${pSp}^{\perp}$, it suffices to show that
$\delta = d_b$ on each ${pap}_{\alpha}$. Each of these is a rank one
operator $\xi_{\alpha} \otimes \eta$ with $\eta = pa\xi_\alpha \in p$,
so we may restrict our attention to such operators. 

Consider, then, $q=\xi_{\alpha} \otimes \eta$, $\eta \in p$. Then $q$ is as in
the statement of the theorem. 

For any $\xi \in p^{\perp}, q \xi=<\xi, \xi_{\alpha}>\eta
$
so
\beqn
\begin{aligned}
b q \xi & ~=~{b_1} q \xi ~=~ <\xi, \xi_{\alpha}>{b_1} \eta  \\
& ~=~ <\xi, \xi_{\alpha}>\delta\left(q_{2}\right) \xi_{0} \\
& ~=~\delta\left(q_{2}\right) q_{1}^{*} q_{\alpha} \xi .
\end{aligned}
\eeqn
Then
\beqn(b q-q b) \xi={b_1} q \xi-q p_{\alpha} b \xi
=\left[\delta\left(q q_{\alpha}^{*} q_{1}\right) q_{1}^{*}
  q_{\alpha}+q q_{\alpha}^{\star} \delta\left(q_{\alpha}\right)-q
  q_{\alpha}^{*} \delta\left(q_{1}\right) q_{1}^{*} q_{\alpha}\right]
\xi 
\eeqn
However condition (\ref{theorem13eq}) says that this is just
$\delta(q) \xi$ which is 
what we require.  

If instead $\xi \in p$ is a unit vector, put $q^{\prime}=\xi_{0}
\otimes \xi \in S$. Then $q q^{\prime}=0$ and $q^{\prime} \xi_{0}=\xi$,
and by (\ref{implement1}) ${b_1} \xi=\delta(q^{\prime}) \xi_{0}$. Hence
\beqn
\delta(q) \xi=\delta(q) q^{\prime} \xi_{0}=-q \delta(q^{\prime})
\xi_{0}=-q {b_1} \xi=-q b \xi=(b q-q b) \xi ,
\eeqn
since $q \xi=0$ and $b p={b_1} p$. Thus $\delta(q)=d_b(q)$ on all of $H$.

We have thus shown that if condition (\ref{theorem13eq}) holds then
$\delta=d_{b} \mid 
p S p^\perp$ and Theorem~\ref{theorem12} shows that $b$ implements $\delta$ on the
rest of $S$, and so the proof is complete.     \qed

\medbreak
\noindent{\bf Remarks.}
\begin{enumerate}[label=\arabic*)]

\item The condition (\ref{theorem13eq}) can be shown to be a special case of
  the product rule (Property~$D_1$) of \cite{spivack1987certain}.
For $\delta$ strongly continuous on $pSp^{\perp}$, Theorem
\ref{theorem13} shows that (\ref{theorem13eq}), Property~$D_1$, and
implementation of $\delta$ are equivalent. 

\item The operator $b$ can be constructed in an 
  identical way for any non-trivial nest algebra $N$ since
  the construction uses only operators in ${pB}({H}) {p}^{\perp}$ and
  these lie in $N$ whenever $p\in \operatorname{lat}(N)$.
Strong continuity and (\ref{theorem13eq}) hold on $N$ because every
derivation of $N$ is implemented
\cite{christensen1,lance,spivack1985derivations}.
  Thus the construction yields an explicit implementing operator for nest
  algebras without distinguishing lattices discrete and continuous at 1.

  \item The definitions and results here
    lend themselves naturally to Banach space generalization if we
    should care to imitate the process used in
    \cite{spivack1982thesis,spivack1985derivations} for nest algebras.  
\end{enumerate}

\section*{Acknowledgments}

I remain deeply grateful to Professor Chris Lance for his kind generosity
as well as to Professor Roger Plymen for their guidance during
the writing of \cite{spivack1982thesis} from which this note is adapted.

\bibliographystyle{plain}
\bibliography{refs}

\end{document}